\documentclass[12pt, twoside, myheadings]{amsart}
\usepackage {amsmath, amscd}
\usepackage {amssymb}
\usepackage {epsfig}
\usepackage {bm}
\usepackage {indentfirst} 

\usepackage{mathtools}

\usepackage{tkz-tab}

\usetikzlibrary{patterns}

\usepackage[all]{xy}

\usepackage{cancel}
\usepackage{ulem}
\usepackage{xcolor}

\usepackage{lipsum}
\usepackage{mwe}
\usepackage{subcaption}

\usepackage[active]{srcltx}

\advance\hoffset by -1.5cm

\numberwithin{equation}{section}

\def\<{\langle}
\def\>{\rangle}

\def\N{{\bf N}}

\newtheorem{lemma}{Lemma}[section]
\newtheorem{theorem}[lemma]{Theorem}
\newtheorem{corollary}[lemma]{Corollary}
\newtheorem{definition}[lemma]{Definition}

\theoremstyle{definition}
\newtheorem{remark}[lemma]{Remark}

\newcommand{\pfend}{\hfill $\Box$ \medskip}

\usepackage{fancyhdr}
\usepackage{graphicx}
\graphicspath{ {./images/} }

\title{Signed representing measures \\ (Berger-type charges) in subnormality \\ and related properties of weighted shifts}
\author{Chafiq Benhida}
\address{UFR de Math\'{e}matiques, Universit\'{e} des Sciences et
Technologies de Lille, F-59655 \newline Villeneuve-d'Ascq Cedex, France}
\email{chafiq.benhida@univ-lille.fr}

\author{Ra\'ul E. Curto}
\address{Department of Mathematics, University of Iowa, Iowa City, Iowa 52242-1419, USA} 
\email{raul-curto@uiowa.edu}

\author{George R. Exner}
\address{Department of Mathematics, Bucknell University, Lewisburg, Pennsylvania 17837, USA}
\email{exner@bucknell.edu}



\begin{document}


\begin{abstract}
In the study of the geometrically regular weighted shifts (GRWS) -- see \cite{BCE3} -- signed power representing measures (which we call Berger-type charges) played an important role. \  Motivated by their utility in that context, we establish a general theory for Berger-type charges. \  We give the first result of which we are aware showing that $k$--hyponormality alone (as opposed to subnormality) yields measure/charge-related information. \ More precisely, for signed countably atomic measures with a decreasing sequence of atoms we prove that $k$-hyponormality of the associated shift forces positivity of the densities of the largest $k+1$ atoms. \ Further, for certain completely hyperexpansive weighed shifts, we exhibit a Berger-type charge representation, in contrast (but related) to the classical L\'{e}vy-Khinchin representation. \  We use Berger-type charges to investigate when a non-subnormal GRWS weighted shift may be scaled to become conditionally positive definite, and close with an example indicating a distinction between the study of moment sequences and the study of weighted shifts.
\end{abstract}

\maketitle

\tableofcontents

\setcounter{tocdepth}{2}


\section{\bf Introduction and  preliminaries}  \label{se:Intro}

Let $\mathcal{H}$ be a separable infinite dimensional complex Hilbert space and $\mathcal{L}(\mathcal{H})$ the collection of bounded linear operators on $\mathcal{H}$. \  Properties such as subnormality, $k$-hyponormality, complete hyperexpansivity, and complete positive definiteness of operators (all definitions reviewed below) have received considerable study. \  A standard starting place for such studies, or for the development of new tools and definitions, has been the class of weighted shifts. \  In \cite{BCE3} the authors studied a new class of shifts -- the geometrically regular weighted shifts (GRWS) -- which provide, by variation in a certain parameter space, examples in these classes. \ Important in that study were signed representing measures (henceforth, Berger-type charges). \  In this paper we return to such GRWS to emphasize the utility of Berger-type charges in that setting and prove as well some more general results suggested from that study.

The organization of this paper is as follows. \  In the remainder of this section we give necessary definitions and background. \  In Section \ref{se:Mainresults} we state our main results. \  In Section \ref{se:Proofanddiscussion} we gives proofs and discussion of these and related results.

We employ $\mathbb{Z}$, $\mathbb{Z}_+$, $\mathbb{R}$, and $\mathbb{R}_+$ for the integers, non-negative integers, real numbers, and non-negative real numbers respectively. \  We define an operator $\Delta$ on sequences as follows:  given a sequence $c = (c_n)_{n=0}^\infty$,  set $\Delta(c)$ by $(\Delta(c))_n = c_{n+1} - c_n$ for all $n$;  we define $\Delta^n$ by $\Delta^0(c) = c$ and $\Delta^{n+1}(c) = \Delta(\Delta^n(c))$. \  The operator $\nabla$ is simply $-\Delta$. \ Recall that a sequence $c$ is said to be completely monotone if, for all $n$ and $k$, $(\nabla^k(c))_n \geq 0$. \ We say $c$ is completely alternating if for all $n$ and $k$, $(\Delta^k(c))_n \geq 0$.

Consider $\ell^2$ with its standard basis $\{e_j\}_{j=0}^\infty$ (note that we begin indexing at zero). \ Given a weight sequence $\alpha: \alpha_0, \alpha_1, \alpha_2,
\ldots$ (which, for almost all questions of interest, we may and do assume is strictly positive), we define the weighted shift $W_\alpha$ on $\ell^2$ by $W_\alpha e_j = \alpha_j e_{j+1}$, and extend
by linearity. \  The \underline{moments} of the shift are $\gamma_0 = 1$ and $\gamma_j =\prod_{i=0}^{j-1} \alpha_i^2$, $j \geq 1$. \  An operator $T$ is contractive if, for all $x$, $\|Tx\| \leq \|x\|$;  for shifts, this is easily that $\alpha_i \leq 1$ for all $i$.

Standard areas of study, often using weighted shifts, have been subnormality and related properties. \    An operator $T$ is normal if $T^* T = T T^*$ and subnormal if it is the restriction of a normal operator to a (closed, linear) invariant subspace. \  One route to showing an operator is subnormal uses the Bram-Halmos characterization employing $k$-hyponormality, $k = 1, 2, \ldots$ (see \cite{Br}). \  The operator matrix positivity condition for $k$-hyponormality simplifies considerably for weighted shifts:  a weighted shift is $k$-hyponormal if and only if the $(k+1)\times (k+1)$ Hankel matrices of moments
\begin{equation}  \label{eq:Hankelmomentmatrix}
\left( \begin{array}{ccccc}
    \gamma_m & \gamma_{m+1} & \gamma_{m+2} & \ldots & \gamma_{m+k}  \\
    \gamma_{m+1} & \gamma_{m+2} & & \ldots & \gamma_{m+k+1} \\
  \gamma_{m+2} & \ldots & & \ldots & \gamma_{m+k+2} \\
     \vdots &  & \vdots  &  & \vdots \\
     \gamma_{m+k} & \gamma_{m+k+1}&  & \ldots & \gamma_{m+2k}
                    \end{array} \right) \geq 0, \quad (m = 0, 1, 2, \ldots)
\end{equation}
are all positive (semi)-definite \cite[Definition 3]{Cu}. \ Another route to subnormality for weighted shifts is the presence of a Berger measure:  the weighted shift $W$ is subnormal if its moment sequence $\gamma$ has a representing (Berger) probability measure $\mu$ (\cite{GW} and \cite[III.8.16]{Con}) compactly supported in $\mathbb{R}_+$ (if $W$ is a contraction, in $[0,1]$) such that
\begin{equation} \label{eq 12}
\gamma_n = \int t^n \, d \mu(t), \quad (n = 0, 1, 2, \ldots).
\end{equation}

This paper concerns itself primarily with the replacement of the (positive) power representing Berger measure, as in (\ref{eq 12}), with a signed power representing measure (hereafter we call this a Berger-type charge). \ We will frequently consider atomic measures, so let $\delta_a$ be the Dirac measure at the point $a$.

We pause for a moment to consider the uniqueness of Berger-type charges. \ Of course it is well-known that moment sequences may be indeterminate (have even positive measure representations that are not unique -- see \cite{Si}, where the key question is the rate of growth of the sequence);  it is known as well that any sequence has a signed representing measure (\cite{Bo}) and the construction used clearly allows for non-uniqueness. \  In the cases we consider, however, the charges are compactly supported, and splitting the charge on its positive and negative sets, and using the uniqueness of (positive) representing measures with compact support (which then satisfy the growth condition) one can show that if we normalize the Berger-type charge so that $\gamma_0 = 1$ as is traditional we indeed have uniqueness. \  (See Remark \ref{remark:signedrepresentations} on page \pageref{remark:signedrepresentations} for a more nuanced discussion of normalization.)

A weighted shift $W_\alpha$ is moment infinitely divisible ($\mathcal{MID}$) -- a property stronger than subnormality -- if  $W_{\alpha^{(p)}}$ is subnormal for every $p > 1$, where $\alpha^{(p)}$ denotes the weight sequence $(\alpha_n^p)_{n=0}^\infty$. \  A study of these shifts, which have extremely nice properties, was initiated in \cite{BCE1} and \cite{BCE2}. \  Recall that a shift is completely hyperexpansive (CHE) if its moment sequence $\gamma$ is completely alternating. \  The moment sequence of such a shift has the classical L\'{e}vy-Khinchin representation (see \eqref{eq:LKrepresentation}). \  A shift is conditionally positive definite (CPD) if positivity of the matrices in \eqref{eq:Hankelmomentmatrix} is replaced by positivity of the matrices, as a quadratic form, but restricting to vectors $v = (v_0, \ldots, v_k)^T$ such that $\sum_{i=0}^k v_i = 0$. \  Equivalently, in the presence of a growth condition satisfied in our setting, $\Delta^2(\gamma)$ is positive definite (see \cite[Prop. 2.2.9]{JJS}). \ Recall that conditional positive definiteness for matrices has been well-studied (see, for example \cite{Bh}, \cite{HJ1}, and \cite{HJ2}). \  In \cite{JJS} and \cite{JJLS} the authors consider in depth the notion of CPD for sequences, for operators in general, and for weighted shifts in particular.

\begin{figure}


{\begin{tikzpicture}[scale=4]

\draw[->] (-1.3,0) -- (1.3,0) ;

\draw[->] (0,-1.3) -- (0,1.3) ;

\draw  [-, dashed] (-1,-1)--(1,-1) --(1,1)--(-1,1)--cycle  ;

\draw (1,0) node[below right] {1} ;

\draw (-1,0) node[below left] {-1} ;

\draw (0,1) node[above left ] {1} ;

\draw (0,-1) node[below left] {-1} ;


\draw [][  domain=-1:1] plot(\x,{\x })node[right, scale=0.8] {$D=N$} ;

\draw [][  domain=-1:1] plot(\x,{-\x })node[right, scale=0.8] {$D=-N$} ;

\draw [blue ][  domain=-2/3:2/3] plot(\x,{3/2*\x })node[above right, scale=0.65] {$D=pN$} ;

\draw [-,blue, dashed][  domain=-0:4/9] plot(\x,{9/4*\x })node[above , scale=0.6] {$D=p^2N$} ;

\draw [-, blue, dashed][  domain=0:8/27] plot(\x,{27/8*\x })node[above, scale=0.8] {} ;

\draw [-, blue, dashed][  domain=0:16/81] plot(\x,{81/16*\x })node[above, scale=0.8] {} ;

\fill[color=red!60] (0,0) -- (-1,-1) -- (-2/3,-1)  -- cycle ;

\fill[  color=red!10 ] (0,0) -- (-2/3,-1)--(-4/9,-1)  -- cycle ;

\fill[color=green] (0,0) -- (-1,-1) -- (-1,0)  -- cycle ;

\fill[green!40] (0,0) -- (-1,0) -- (-1,1)  -- cycle ;

\fill[blue!70] (0,0) -- (0,1) -- (-1,1)  -- cycle ;

\draw [green] (0,0) -- (1,1);

\draw (-0.6, -0.3) node [ ] {I};

\draw (-0.6, 0.3) node [ ] {II};

\draw (-0.3, 0.6) node [ ] {III};

\draw (0.3, 0.6) node [ ] {IV};

\draw (0.6, 0.3) node [ ] {V};

\draw (0.6, -0.3) node [ ] {VI};

\draw (0.3, -0.6) node [ ] {VII};

\draw (-0.3, -0.6) node [ ] {VIII};

\draw (-0.75, -0.9) node [scale=0.5 ] {VIII.A};

\draw (-0.50, -0.9) node [scale=0.5 ] {VIII.B};


\end{tikzpicture}}

\caption{Magic Square}
\label{MagicSquare}
\end{figure}
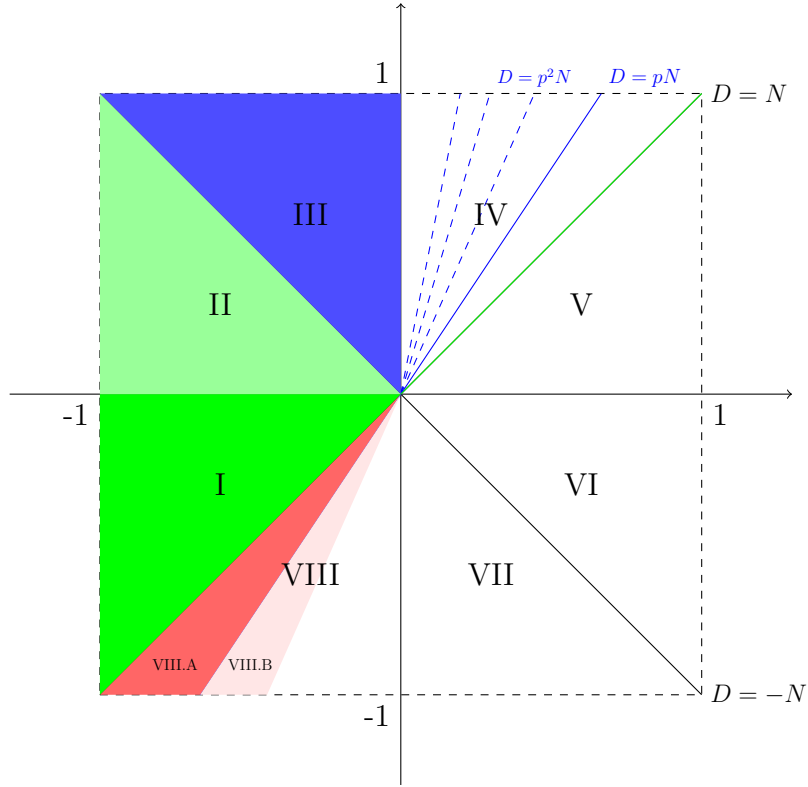

In \cite{BCE3} the authors study a class of weighted shifts we call geometrically regular weighted shifts (GRWS). \  Given $p > 0$ and $(N, D)$ in the open unit square $(-1,1)\times (-1,1)$ (which we have come to call the Magic Square), the shift has weights
$$\alpha(N,D) := \sqrt{\frac{p^n + N}{p^n + D}}.$$
Subsectors of the Magic Square (Figure \ref{MagicSquare}) yield examples in various classes of interest.
\begin{itemize}
\item Sector I:  $\mathcal{MID}$, with the stronger property that the weights squared are interpolated by a Bernstein function.

\item Sector II:  $\mathcal{MID}$, with the stronger property that the weights squared are interpolated by a log Bernstein function.

\item Sector III:  Subnormal

\item Sector IV:  Subnormal (finitely atomic) on the special lines $D = p^n N$;  varying $k$- hyponormal but not $(k+1)$-hyponormal in subsectors, countably atomic

\item Sectors V-VII:  \textit{terra incognita}

\item Sector VIIIA:  Completely hyperexpansive

\item Sector VIIIB:  See Subsection \ref{BergerchargesandCPD}

\end{itemize}

Every moment sequence from a GRWS has a Berger-type charge arising as follows (\cite{BCE3}):  first, define multipliers $m_i$ by $m_0 := 1$ and
$$m_i := \frac{p(D - p^{i-1}N)}{p^i - 1}, \quad (i = 0, 1, 2, \ldots).$$
Define $c_n$ for $n = 0, 1, 2, \ldots$ by
$$c_n := \prod_{i=0}^n m_i,$$
and set
$$a(N,D) := \frac{1}{\sum_{n=0}^\infty c_n}.$$
Then the (signed) measure $\mu$ defined by
$$\mu := a(N,D)\left(\sum_{i=0}^\infty c_i \delta_{\frac{1}{p^i}}\right)$$
is a representing Berger-type charge for $\gamma$. \ Further, it satisfies $\int 1 \, d \mu = 1$, and thus is a probability measure and thus a Berger measure if it is a measure.

For those values of the parameters not yielding subnormality, there are negative densities for the Berger-type charge. \  The patterns of these densities are as follows, taking things from right to left (the decreasing order of the atoms $1$, $\frac{1}{p}$, $\frac{1}{p^2}$, \ldots) :
\begin{itemize}  \label{listofdensities}

\item Sectors I, II, III:  $+, +, +, \ldots$;

\item Sector IV:  On special lines, $+, +, \ldots, +, 0, 0, 0, \ldots$, because some multiplier $m_i = \frac{p(D - p^{i-1}N)}{p^i - 1}$ becomes zero;  in the subsectors, $+, +, \ldots, +, -, +, -, +, \ldots$

\item Sectors V-VII:  mixed, but with both several positive and negative densities

\item Sector VIIIA:  $+, -, -, -, \ldots$

\item Sector VIIIB:  $+, -, +, +, \ldots$

\end{itemize}

\section{\bf Main results}  \label{se:Mainresults}

\begin{theorem}
Let $\sigma$ be a finite positive measure on $[0,1)$ and let $C \in \mathbb{R}_+$ be such that $C - \sigma([0,1)) = 1$. \  Let $\mu$ be the signed measure defined by $\mu = C \delta_1 - \sigma$. \  If $\gamma = (\gamma_n)_{n=0}^\infty$ is the sequence of moments of $\mu$, and $W_\alpha$ is the weighted shift with weight sequence $\alpha$ induced by the moments $\gamma_n$, so  $\alpha_n = \sqrt{\frac{\gamma_{n+1}}{\gamma_n}}$ for all $n$, then $W_\alpha$ is completely hyperexpansive. \  Further, in the L\'{e}vy-Khinchin representation of $\gamma$,
$$
\gamma_n = 1 + b n + \int_{[0,1)} (1-t^n) \, d \nu(t) \quad (n = 0, 1, 2, \ldots),
$$
we have $\nu = \sigma$.
\end{theorem}

\begin{theorem}
Let $\gamma$ be the moment sequence of a CHE weighted shift, and let $\mu$ be the Berger measure of the associated subnormal weighted shift $W$ with moment sequence $\Delta \gamma$. \ Assume now that $W$ is contractive. \ Then $\gamma$ admits a Berger-type charge representation if and only if 
$$
\int_0^1 \dfrac{1}{1-x} d\mu(x) < \infty.
$$
\end{theorem}

\begin{theorem}
Let $\mu$ be a (potentially signed) measure of the form $\sum_{i=1}^\infty a_i \delta_{r_i}$ where $a = (a_i)_{i=1}^\infty \in \ell^1$, $a_1 > 0$, and the $r_i$ form a strictly decreasing positive sequence. \  Suppose as well that $\sum_{i=1}^\infty a_i = 1$. \  Let $\gamma = (\gamma_n)_{n=0}^\infty$ be the moment sequence for $\mu$, so $\gamma_n = \int_{\mathbb{R}_+} t^n \, d \mu(t)$. \   Let $W_\alpha$ be the weighted shift with weight sequence $\alpha$ induced by the moments $\gamma_n$. \  Then if the restriction of $W_\alpha$ to some canonical invariant subspace $\bigvee \{e_i: i \geq I \}$ is $k$-hyponormal for some $k$, the first $k+1$ non-zero $a_i$ are strictly positive.
\end{theorem}

This is the only result of which we are aware that deduces, from the $k$--hyponormality condition rather than subnormality, information about a measure. \ We remark that we obtain the result via a more technical approach which yields in fact further information (see Lemma \ref{le:workingdeterminants} and Corollary \ref{cor:densityinfonokHN}). \  Note also that $a_1 < 0$ does not allow for an all-positive moment sequence.

\begin{theorem}  \label{th:CPDfirststatement}
Let $\mu$ be a countably atomic (potentially signed) measure with support in $[0, \infty)$. \  If $\mu$ has all negative or all positive densities, the shift $W_\alpha$ is subnormal and hence conditionally positive definite (CPD). \  If $\mu$ has two atoms with one density of each sign, then there are exactly two (non-zero) constants which, multiplying the weights of $W_{\alpha}$, make it CPD. \  If $\mu$ has more than two atoms, with either all positive densities but one or all negative densities but one, there exists exactly one (non-zero) multiplier of the weights of $W_\alpha$ producing a CPD shift. \  If $\mu$ has at least two atoms with positive densities and at least two with negative densities, no non-zero multiple of $W_\alpha$ is CPD.
\end{theorem}

We also discuss a ``CPD-like'' representation for completely hyperexpansive weighted shifts which relies upon a signed measure (as with trivial exceptions CHE weighted shifts are not CPD).

\section{\bf Proofs and discussion}  \label{se:Proofanddiscussion}

\subsection{A Berger-type charge for some CHE shifts}

In \cite{BCE3} were considered weighted shifts with weights
$$\alpha_n = \sqrt{\frac{p^n + N}{p^n + D}}$$
where $p > 1$ and $(N,D)$ lies in the ``Magic Square'' $(-1,1)\times(-1,1)$. \ If $(N,D)$ falls in a certain sector (Sector VIIIA) of parameter space, the shifts produced are completely hyperexpansive (CHE) shifts. \ In the generic case (setting aside some special values of the parameters) the Berger-type charge for the relevant shifts has an atom with positive density at $1$, and all the densities for the atoms at $\frac{1}{p^i}$, $i = 1, 2, \ldots$,  are strictly negative. \ This suggests a generalization to a subclass of the class of CHE shifts.

Recall that the moment sequence $\gamma$ of a completely hyperexpansive weighted shift has, since that moment sequence is completely alternating, the L\'{e}vy-Khinchin representation (see, for example, \cite[Ch. 4, Prop. 6.11]{BCR})
\begin{equation}  \label{eq:LKrepresentation}
\gamma_n = a + b n + \int_{[0,1)} (1-t^n) \, d \nu(t) \quad (n = 0, 1, 2, \ldots),
\end{equation}
where $a = \gamma_0 = 1$, $b \geq 0$, and $\nu$ is a positive Radon measure on $[0,1)$.

\begin{theorem} \label{thm31}
Let $\sigma$ be a finite positive Radon measure on $[0,1)$ and let $C:= 1 + \sigma([0,1))$. \  Let $\mu$ be the signed measure defined by $\mu = C \delta_1 - \sigma$. \  If $\gamma = (\gamma_n)_{n=0}^\infty$ is the sequence of moments of $\mu$, and $W_\alpha$ is the weighted shift with weight sequence $\alpha$ induced by the moments $\gamma_n$, so  $\alpha_n = \sqrt{\frac{\gamma_{n+1}}{\gamma_n}}$ for all $n$, then $W_\alpha$ is completely hyperexpansive. \ Further, in the L\'{e}vy-Khinchin representation of $\gamma$, we have $\nu = \sigma$, $a = 1$, and $b = 0$.
\end{theorem}

\noindent Proof. \ It is well-known that a sequence $s = (s_n)_{n=0}^\infty$ is completely alternating if and only if the sequence $\Delta(s)$ is completely monotone (\cite[Ch. 4, Lemma 6.3]{BCR}). \  But it is easy to see that computing $\Delta$, in measure terms, for the sequence of moments of $\mu$ is to move to the measure $\mu'$ where $\mu'(t) =-(t-1)\sigma(t)$. \  Since this is a (finite) positive measure, and with closed support contained in $[0,1]$, the resulting weighted shift is contractive and its moment sequence is completely monotone (\cite{Ag}), as desired. \ The equality $\nu = \sigma$, $a = 1$, and $b = 0$ follows from the computation (from $\mu'$ to $\nu$ in \cite[page 3747]{At}).  \pfend

That the Berger-type charge representation forces $b=0$ in the L\'{e}vy-Khinchin representation (\ref{eq:LKrepresentation}) indicates that this representation does not capture all CHE shifts, as is obvious since the well-known Dirichlet shift arises from $a = 1$, $b=1$, and $\nu = 0$. \ In the sequel, we present a result that exhibits a class of CHE weighted shifts that do admit a representation of the form (\ref{eq:LKrepresentation}).

Consider a CHE weighted shift and let $\gamma$ be its moment sequence. \ Assume that $\gamma$ has a Berger-type charge representation. \ Thus, there exist two positive Borel measures $\nu_1$ and $\nu_2$ supported on, say an interval $[0,a]$, with $a\geq 1$, such that

$$
\gamma_n= \int_0^a x^n d\nu_1(x)-\int_0^a x^n d\nu_2 (x) \quad (n = 0, 1, 2, \ldots).
$$

(As we will show below, without loss of generality we may take $a=1$.)

Since $\gamma$ is completely alternating, $\Delta\gamma$ is (modulo normalization) the moment sequence of a contractive subnormal weighted shift. \ Then there exists a positive Borel measure $\mu$ supported on $[0,1]$ and a positive constant c such that

$$
\Delta\gamma(n)=c \int x^n\chi_{[0,1]}d\mu(x) \quad (n = 0, 1, 2, \ldots).
$$

This implies

$$
c \int x^n\chi_{[0,1]}d\mu(x)=\int_0^a x^n(x-1) d\nu_1(x)-\int_0^a x^n(x-1) d\nu_2(x) \quad (n = 0, 1, 2, \ldots).
$$

From this, we can now deduce that 

$$
\int_{(1,a]} x^nd\nu_1 (x)= \int_{(1,a]}  x^nd\nu_2 (x) \quad (n = 0, 1, 2, \ldots),
$$
and consequently, as announced above, we may always assume that $a=1$.

Now,
$$
c \int_0^1 x^nd\mu(x)=
\int_0^1 x^n(1-x) d\nu_2(x) -\int_0^1 x^n(1-x) d\nu_1(x) \quad (n = 0, 1, 2, \ldots),
$$
or equivalently,
$$
\int_0^1 x^n[cd\mu(x)+ (1-x) d\nu_1(x)] = \int_0^1 x^n(1-x) d\nu_2(x) \quad (n = 0, 1, 2, \ldots),
$$
implying that
$$
c \cdot d\mu(x)+ (1-x) d\nu_1(x)=(1-x)d\nu_2(x).
$$
From this, we readily infer that 
\begin{equation}\label{necessary}
\int_0^1\frac{1}{1-x} d\mu(x)<\infty.
\end{equation}
As a result, we obtain (\ref{necessary}) as a necessary condition for the existence of a Berger-type charge; in other words, the Berger measure of the contractive subnormal weighted shift associated with the sequence $\Delta \gamma$ must satisfy (\ref{necessary}). \ We will now prove the converse, i.e., (\ref{necessary}) is also a sufficient condition for $\gamma$ to admit a Berger-type charge representation of the form (\ref{eq:LKrepresentation}).  

\begin{theorem} \label{sufficient} \ Given $\mu$ satisfying (\ref{necessary}), the moment sequence $\gamma$ given by
$$
\gamma_n = 1+  c\int_0^1 (1+x+\dots +x^{n-1})d\mu(x) \quad (n = 0, 1, 2, \ldots)
$$
has a Berger-type charge representation.
\end{theorem}

\noindent Proof. \ Recall that
$$
\begin{array}{lcl}
 \gamma_n &=& 1+  c\int_0^1 (1+x+\dots +x^{n-1})d\mu(x)\\ & & \\
 &=& 1+ c\int_0^1 \frac{1-x^n}{1-x}d\mu(x)\\ & & \\
   &=& \Big[1+ c \int_0^1\frac{1}{1-x} d\mu(x) \Big] - \int_0^1 x^n c \cdot \dfrac{d\mu(x)}{1-x} \quad (n = 0, 1, 2, \ldots)
 \end{array}
$$
Therefore, the weighted shift whose moment sequence is $\gamma$ has, as a representing measure, the Berger-type charge
$$
\Big[1+ c\int_0^1\frac{1}{1-x} d\mu(x) \Big]\cdot \delta_1- c \cdot \chi_{[0,1]} \dfrac{d\mu(x)}{1-x}.
$$ \pfend 

As an application, we see that none of the weighted shifts $W_{(\sqrt{\frac{n+m}{n+1}})_{n\in\N}}$ (where $1<m\leq 2$)
 has a Berger-type charge representation.


Even with the limitation that such a CHE representation $\mu=C \delta_1 - sigma$, as in Theorem \ref{thm31}, is not available for all CHE shifts, the use of such a representation facilitates the taking of Schur products of shifts (convolutions of their measures) in the case of the GRWS and being able to ``spot'' when the result is, or is not, CHE. \  We turn next to a brief discussion of this use of measures to ``detect'' or reject membership in desirable classes.

\subsection{Berger-type charges in detection/rejection}

A Berger-type charge carries information about membership, or failure of membership, in various classes. \  For example, suppose one has a Berger-type charge $\mu$ which is finite or countably infinite and with moment sequence $\gamma$. \ If the support is contained in $[0,1]$ but there is at least one negative density, the computation above shows that the measure for $\Delta(\gamma)$ is unlikely to be a positive measure since that negative density is preserved. \  (See Subsection \ref{BergerchargesandCPD} for a more nuanced discussion of this approach.)  If the support is not contained in $[0,1]$, it may be relatively easy to show that $\Delta(\gamma)$ is not completely monotone since its measure computation is easy.

As another example, if we take the Schur product of two shifts (the convolution of their measures, which is particularly easy in the atomic case since $\delta_a \ast \delta_b$ is just $\delta_{ a b}$) and obtain a measure of the sort in the theorem above, we know we have obtained a CHE shift. \  Should we obtain a measure with all positive densities, we have obtained a subnormal shift.

\begin{figure}[ht]


{\begin{tikzpicture}[scale=3]

\draw[->] (-1.5,0) -- (1.5,0) ;

\draw[->] (0,-1.5) -- (0,1.5) ;

\draw  [-, dashed] (-1,-1)--(1,-1) --(1,1)--(-1,1)--cycle  ;

\draw [-, dashed , gray][  domain=1:-1] plot(\x,{-\x })node[left, scale=0.8] {$y=-x$} ;

\draw [fill=blue!40]  [-,  blue!40] (-1,-1)-- (-1/4, -1/4)--(-1/4, 0)--(0, 1/4)--
(1,1/4)--(1,0)--(1, -1)--cycle  ;

\draw (-1/4, 1/2) node[blue, scale=0.5] {$\bullet$} ;

\draw (-1/4,1/2) node[ left, scale=0.5] {$(N,D)$} ;

\draw (-1/4,0) node[ below left, gray, scale=0.5] {$N$} ;

\draw (0,1/2) node[ above left, gray, scale=0.5] {$D$} ;

\draw [-, thick , blue!40]  (-1/4, 1/2)--  (1, 1/2);

\draw [-, thick , blue!40]  (-1/4, 1/3)--  (1, 1/3);

\draw [-, thick , blue!40]  (-1/4, 2/9)--  (1, 2/9);

\draw [-, thick , blue!40]  (-1/4, 4/27)--  (1, 4/27);

\draw [-, thick , blue!40]  (-1/4, 8/81)--  (1, 8/81);

\draw [-, thick , blue!40]  (-1/4, 16/243)--  (1, 16/243);

\draw [-, thick , blue!40]  (-1/4, 32/729)--  (1, 32/729);

\draw [fill=blue!40]  [-,  blue!40]
 (0,1/4)--
(-2/9, 1/4)-- (-1/4, 2/9)--(-1/36, 0)
--cycle  ;

\draw [fill=blue!40]  [-,  blue!40]
 (0,1/4)--
(-4/27, 1/4)-- (-1/4, 4/27)--(-11/108, 0)
--cycle  ;

\draw [fill=blue!40]  [-,  blue!40]
 (0,1/4)--
(-8/81, 1/4)-- (-1/4, 8/81)--(-49/324, 0)
--cycle  ;

\draw [fill=blue!40]  [-,  blue!40]
 (0,1/4)--
(-16/243, 1/4)-- (-1/4, 16/243)--(-179/972, 0)
--cycle  ;

\draw [fill=blue!40]  [-,  blue!40]
 (0,1/4)--
(-32/729, 1/4)-- (-1/4, 32/729)--(-601/2954, 0)
--cycle  ;

\draw [fill=blue!40]  [-,  blue!40]
 (0,1/4)--
(-64/2187, 1/4)-- (-1/4, 64/2187)--(-2123/8748, 0)
--cycle  ;


\draw [-, thick, blue!40][  domain=0:1] plot(\x,{\x })node[right , gray, scale=0.8  ] {$y=x$} ;

\draw [-, thick, blue!40][  domain=0:1] plot(\x,{2/3*\x })node[right , gray, scale=0.6  ] {$y=\frac{1}{p}x$} ;

\draw [-, thick, blue!40][  domain=0:1] plot(\x,{4/9*\x })node[right, gray, scale=0.6   ] {$y=\frac{1}{p^2}x$} ;

\draw [-, thick, blue!40][  domain=0:1] plot(\x,{8/27*\x })node[right, gray, scale=0.5   ] {$y=\frac{1}{p^3}x$} ;

\draw [-, thick, blue!40][  domain=0:1] plot(\x,{16/81*\x })node[right, gray, scale=0.4   ] {$y=\frac{1}{p^4}x$} ;


\end{tikzpicture}}

\caption{$\mathcal{SQ}$: $(N,D)$ in Sector
 III}
\label{SSectorIII}
\end{figure}
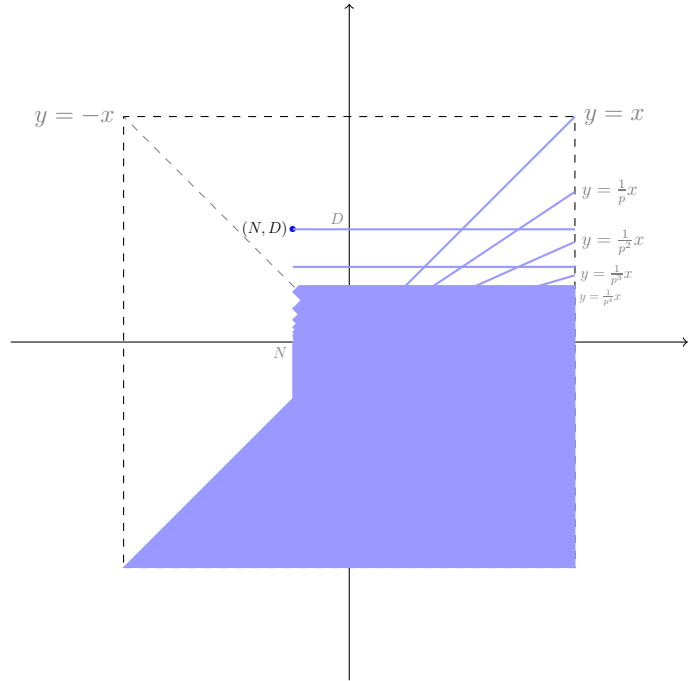

In \cite{BCE4} we were concerned with Schur quotients of the GRWS: starting with a base point $(N,D)$, what are the $(M,P)$ such that the shift with the quotient weights $\frac{\alpha(N,D)}{\alpha(M,P)}$ (which, due to symmetry, may be expressed as a product of such weights), are, or are not, subnormal?  A particularly difficult case is captured in the diagram in Figure \ref{SSectorIII}, where the base point $(N,D)$ is in Sector III (this diagram is Figure 4 in \cite{BCE4}, with the blue region indicating $(M,P)$ yielding subnormality of the Schur quotients).

For points in the interiors of the quadrilaterals and similar regions in the first quadrant, it is natural to wonder whether these also yield safe subnormal quotients or not. \  As noted in \cite{BCE4} one can show, as an example in one such quadrilateral, that not all do. \  However, the natural approach of testing for failure of $k$-hyponormality for perhaps large $k$ becomes computationally intractable even with the aid of \textit{Mathematica} (\cite{Wol}). \  The standard tests for positivity of the needed moment matrices involve the Nested Determinant Test, and the determinants are computationally intensive. \  But the convolution of the measures (using again from symmetry that quotients may be expressed as products) is essentially no worse than the multiplication of generating functions built from the measures. \  As soon as one spots a negative density one has a failure of subnormality (while a more nuanced discussion as in Section \ref{BergerchargesandCPD} is required, the subtleties there don't occur in this case). \  One has found a point not yielding subnormality, and if one wishes to test for failure using $k$-hyponormality one has a useful candidate (and, in fact, the result of the next subsection gives at least an indication of what $k$ is likely to be required).

\subsection{$k$-hyponormality and densities}

In the discussion in \cite{BCE3} of shifts with weights
$$\alpha_n = \sqrt{\frac{p^n + N}{p^n + D}},$$
certain values of the parameters yielded, in subsectors of Sector IV, shifts which were $k$-hyponormal but not $(k+1)$-hyponormal for a desired choice of $k$. \  The Berger-type charges for those shifts are (except for special cases on certain special lines in parameter space) countably atomic and have the property that where there is $k$- but not $(k+1)$-hyponormality the densities for the atoms at $1, \frac{1}{p}, \ldots, \frac{1}{p^k}$ are all positive, the density at $\frac{1}{p^{k+1}}$ is negative, and after that the densities alternate in sign. \  This, as well as striking similarities in certain graphs testing $k$-hyponormality and signs of densities using \textit{Mathematica}(\cite{Wol}) in connection with the Schur quotients in \cite{BCE4}, suggest that this might be generalized.

Note first, however, that it is easy to produce an example of an even finitely atomic measure with geometric support $\{1, \frac{1}{p}, \ldots, \frac{1}{p^N}\}$ with the right number of positive densities as one works right to left but not the right $k$-hyponormality. \  Perhaps unsurprisingly, the route to such an example is a number of positive densities small in magnitude and then a negative density much larger in magnitude. \  A string of positive densities need not produce the $k$-hyponormality most optimistically anticipated on the basis of the results from Sector IV.

However, the result does generalize in the other direction. \  We first need a definition.

\begin{definition}  \label{def:asymposkdets}
Suppose $W_\alpha$ is a weighted shift with positive weights $\alpha$ and moment sequence $\gamma = (\gamma_n)_{n=0}^\infty$. \  For some $k = 1, 2, \ldots$, we say that $W_\alpha$ has asymptotically positive $k$-determinants if, with the usual Hankel moment matrix
$$M_n^{k} :=
\left( \begin{array}{ccccc}
    \gamma_n & \gamma_{n+1} & \gamma_{n+2} & \ldots & \gamma_{n+k-1}  \\
    \gamma_{n+1} & \gamma_{n+2} & & \ldots & \gamma_{n+k} \\
  \gamma_{n+2} & \ldots & & \ldots & \gamma_{n+k+2} \\
     \vdots &  & \vdots  &  & \vdots \\
     \gamma_{n+k-1} & \gamma_{n+k}&  & \ldots & \gamma_{n+2k-2}
                    \end{array} \right),$$
there exists $N$ so that $\det(M_n^{k}) > 0$ for all $n \geq N$.
\end{definition}

Note that these matrices $M_n^k$ are of size $k \times k$, although, as usual, $k$-hyponormality or its potential asymptotic versions considers moment matrices of size $(k+1)\times(k+1)$.

\begin{lemma}  \label{le:workingdeterminants}
Let $\mu$ be a (potentially signed) measure of the form $\sum_{i=1}^\infty a_i \delta_{r_i}$ where $a = (a_i)_{i=1}^\infty \in \ell^1$, $a_1 > 0$, and the $r_i$ form a strictly decreasing positive sequence. \  Suppose as well that $\sum_{i=1}^\infty a_i = 1$. \  Let $\gamma = (\gamma_n)_{n=0}^\infty$ be the moment sequence for $\mu$, so $\gamma_n = \int_{\mathbb{R}_+} t^n \, d \mu(t)$, and let $W_\alpha$ be the weighted shift with weight sequence $\alpha$ induced by the moments $\gamma_n$. \  Then for any $k$, $k = 1, 2, \ldots$, $W_\alpha$ has asymptotically positive $k$-determinants if and only if the first $k$ non-zero $a_i$ are strictly positive.
\end{lemma}

\noindent Proof. \ To ease the notation we simply remove any atom $r_i$ whose density is zero. \  We begin by considering the case in which the measure is finitely atomic with $N$ atoms. \  We seek to show that for some $1 \leq k \leq N$, (strict) positivity of terms in the tail of the sequence
$$ \det (M^{k}_n)$$
forces positivity of the product $a_1 a_2 \cdots a_{k}$.

Recall that we have assumed $a_1 > 0$, which is necessary to have even the possibility of all positive moments;  this is the $k=1$ version of the lemma. \  Note also that for the sizes of matrices we are considering, and for the finitely atomic case, we aren't going to have some asymptotically positive $k$-determinants with $k > N$ since in this case such determinants are zero by the rank/column dependence arguments in \cite{CF}.

Denote the entries of the matrix $M^{k}_n$ by  $M_{i,j}$, temporarily suppressing $k$ and $n$. \  Then of course we know that as usual
\begin{equation}  \label{eq:bigdeterminant}
\det M^{k}_n = \sum_{\pi \in P_{k}} \left(\mbox{\rm sgn}(\pi) \prod_{i=1}^{k} M_{i, \pi(i)}\right)
\end{equation}
where $P_{k}$ is the set of permutations and $\mbox{\rm sgn}$ is the sign of the permutation.
If we think about some such term
\begin{equation}  \label{eq:terminDET}
\prod_{i=1}^{k} M_{i, \pi(i)},
\end{equation}
since each of the $M_{i, \pi(i)}$ arises as a sum like
$$\sum_{j=1}^N a_j {r_j}^m$$
for some $m$ satisfying $n \leq m \leq n + 2k -2$, the product in \eqref{eq:terminDET} is very large. \  If expanded out completely it will be a sum of sums, containing terms which are products of $k$ of the $a_j$ (perhaps with repetitions) and various powers of the $r_j$, corresponding to the $k$ choices of which of the $a_j$-related terms to pick from the entries in the various $M_{i, \pi(i)}$. \  Obviously then the determinant as a whole is a sum of such terms.

Let $\mathcal{C} = \{(j_1, j_2, \ldots, j_m\}$ be some subset of $\{1, 2, \ldots, N\}$, and let $S_\mathcal{C}$ be the collection of all the terms in the expanded determinant in which each, but only, $j$ in $\mathcal{C}$ appears -- sometimes we will be a little careless and think of $\mathcal{C}$ as being the $a_j$ or equivalently $r_j$ instead -- so we might, for example, consider $S_\mathcal{C}$ to be all the terms in the determinant in which all of $a_1$, $a_3$, and $a_4$ and only these occur. \ We have a succession of claims.

\medskip

\noindent Claim 1.  If the cardinality of $\mathcal{C}$ is less than $k$, then the sum of all the terms in $S_\mathcal{C}$ is zero.

\medskip

\noindent Proof of Claim 1.   We use a finite induction. \  Suppose first that $\mathcal{C}$ has but one element $j$. \  Then the terms in $S_{\{j\}}$ are exactly those appearing in  some $k \times k$ Hankel moment matrix determinant (with the appropriate upper-left hand entry $\gamma_n$) from the one-atomic measure in which every $a_i$ with $i \neq j$ is set to zero. \  It is known by a rank (column dependence) argument (see \cite{CF}) that this determinant is zero since $k > 1$, so the sum of the terms in the determinant arising from $S_{\{j\}}$ is zero. \  Suppose now that $C = \{j_1, j_2\}$ with these distinct and with $k \geq 3$. \  The determinant of the $k \times k$ Hankel moment matrix from the two-atomic measure with atoms at only $r_{j_1}$ and $r_{j_2}$ is zero, and hence the sum of its terms is zero. \  But these are exactly the terms in $S_{\{j_1\}} \cup S_{\{j_2\}} \cup S_{\{j_1, j_2\}}$. \  Since their sum is zero, and we know that the sum of the terms from each of the first two sets is zero, we must have that the sum of the terms in $S_{\{j_1, j_2\}}$ is zero. \  We may iterate this argument to yield the result.

Therefore the determinant in \eqref{eq:bigdeterminant}, after removing such zero sums, consists of sums, each associated with a subcollection of the $a_j$, all distinct, of size $k$, of the terms arising from the subcollection.

Fix for the moment some such subcollection $\mathcal{C} = \{j_1, j_2, \ldots, j_{k} \}$.

\medskip

Claim 2.  Each term in the sum associated with this subcollection is of the form
$$a_{j_1}a_{j_2} \cdots a_{j_{k}} \cdots r_{j_1}^n r_{j_2}^n \cdots  r_{j_{k}}^n \cdot K$$
where $K$ is a constant not containing any of the $a_j$, not containing $n$, and consisting of a (signed) product of integer powers of the $r_{j_i}$ with each such power no smaller than zero and no larger than $2k-2$.

\medskip

\noindent Proof of Claim 2.  Each term in this subcollection arises from one of the $\prod_{i=1}^{k} M_{i, \pi(i)}$, as in \eqref{eq:terminDET}, where in addition we have chosen, in picking terms to use, distinct $j_i$ from $\mathcal{C}$ for all the terms. \  Each $a_{j_i}$ will appear, as will each $r_{j_i}$ to some power at least $n$ but no more than $n+2k-2$ because of the size of the matrix. \  If we factor out all the $a_{j_i}$ as well as each $r_{j_i}^n$, what we leave behind is powers of each $r_{j_i}$ to some constant in the claimed range.

Consider now all the terms as in Claim 2 in the portion $S_{\mathcal{C}}$ of the determinant arising from the subcollection $\mathcal{C} = \{j_1, j_2, \ldots, j_{k} \}$.

\medskip

\noindent Claim 3.  The sum of all the terms in $S_{\mathcal{C}}$ contributing to the determinant as just above in Claim 2 is of the form
$$a_{j_1}a_{j_2} \cdots a_{j_{k}} \cdots r_{j_1}^n r_{j_2}^n \cdots  r_{j_{k}}^n \cdot L$$
where $L$ is a positive constant which is the sum of the $K$'s as in Claim 2.

\medskip

\noindent Proof of Claim 3.  If we add the terms, and then factor out the common term
$$a_{j_1}a_{j_2} \cdots a_{j_{k}} \cdot r_{j_1}^n r_{j_2}^n \cdots  r_{j_{k}}^n,$$
 we certainly get a sum for some such $L$. \  This $L$ is positive because suppose we set all $a_j$ not in the collection  $\mathcal{C}$ to zero, so that we are dealing with an $k$-atomic measure. \  Suppose we further temporarily assume that $a_{j_1},a_{j_2}, \ldots a_{j_{k}}$ are all strictly positive. \  Then the sum in Claim 3 is the determinant of the $k\times k$ moment matrix for this measure, with $\gamma_n$ in the upper left, because all we have omitted from this determinant are terms which contain repeated elements from the collection $\mathcal{C}$ (and hence not all its $r_j$ appear) and we know the sum of all these is zero, by Claim 1. \  Since this is then the determinant of an $k \times k$ moment matrix of a shift arising from $k$ atoms (and a positive measure, not a signed measure), the determinant is strictly positive. \  But since the $a_{j_i}$ are assumed positive, and the various $r_{j_i}$ are positive, the term $L$ must be positive in this case. \  But the term $L$ is universal in that it has nothing to do with the assumed signs of the $a_{j_i}$, and therefore it is always positive, as claimed.

 Note for future use that this $L$ is bounded by the sum of the absolute values of all its contributors, which sum is no larger than $k! \prod_{i \in \{j_1, \ldots, j_k\}}r_{i}^{(2k-2)} \leq k! r_1^{k(2k-2)}$. \  Observe this bound is independent of $n$.

 Let $D$ be the determinant of the $k \times k$ Hankel matrix of moments with $a_i$ set to zero for $i = k+1, \ldots, N$;  equivalently, this is the contribution to the total determinant from the sum of the terms in $S_{\{1, 2, \ldots, k\}}$. \  We seek to show this dominates, as $n$ becomes large, the sum of the contributions to the total determinant from all the other $S_\mathcal{C}$. \  Let $B$ denote the maximum $\max_{1 \leq j \leq N} |a_j|$.

 \medskip

\noindent Claim 4.   Consider first the contribution to the total determinant of the terms in some $S_\mathcal{C}$ where $\mathcal{C}$ contains exactly one $j > k$, say $j_1$, and $k-1$ choices of $j$ for $j \leq k$, say $j_2, \ldots, j_k$. \  We claim this contribution is bounded by
\begin{equation}  \label{eq:boundonesmall}
|a_{j_1}|B^{k-1}r_1^n r_2^n \cdots r_{k-1}^n r_{j_1}^n k! r_1^{(2k-2)}.
\end{equation}

\noindent Proof of Claim 4.  By Claim 3 this contribution to the determinant has the form
$$a_{j_1}a_{j_2} \cdots a_{j_{\hat{N}}} \cdots r_{j_1}^n r_{j_2}^n \cdots  r_{j_{k}}^n \cdot L.$$
Note that this is bounded above by
\begin{equation}
a_{j_1}a_{j_2} \cdots a_{j_{\hat{N}}} \cdots r_{j_1}^n r_{j_2}^n \cdots  r_{j_{k}}^n \cdot L \leq |a_{j_1}|B^{k-1}r_1^n r_2^n \cdots r_{k-1}^n r_{j_1}^n k! r_1^{k(2k-2)},
\end{equation}
where we have used the definition of $B$, the fact that the $r_i$ are decreasing, and the bound given for $L$ after the proof of Claim 3.

It is clear that $D$ dominates this expression for $n$ large since $r_{j_1} < r_k$. \  Further, since we are considering now the case in which the number of atoms is finite, $D$ dominates the sum of all the contributions to the determinant of this form as $n$ becomes large, since there are but finitely many.

\medskip

\noindent Claim 5.  The contribution to the total determinant of the terms in some $S_\mathcal{C}$ where $\mathcal{C}$ contains at least two $j$ larger than $k+1$ -- let $j_1$ be the least of these and $j_2$ the next largest -- and let $j_3, \ldots, j_k$ be the remaining $j$ in $\mathcal{C}$, is bounded by
\begin{equation}  \label{eq:boundfortwosmall}
|a_{j_1}| |a_{j_2}|B^{k-2}r_1^n r_2^n \cdots r_{k-2}^n r_{j_1}^n r_{j_2}^n k! r_1^{k(2k-2)}.
\end{equation}

\noindent Proof of Claim 5.  Consider the contribution as above. \   A repetition of the argument in the proof of Claim 4 above shows that this contribution to the determinant is bounded above by
$$|a_{j_1}| |a_{j_2}|B^{k-2}r_1^n r_2^n \cdots r_{k-2}^n r_{j_1}^n r_{j_2}^n k! r_1^{k(2k-2)}.$$

Again, since the number of atoms is finite, there are only finitely many such contributions, and thus $D$ dominates the sum of all these contributions to the determinant as $n$ becomes large.

Since as noted above other contributions to the determinant arising from $\mathcal{C}$ with fewer than $k$ elements are zero, $D = a_{1}a_{2} \cdots a_{k} \cdots r_{1}^n r_{2}^n \cdots  r_{k}^n \cdot L$, with $L$ positive, dominates all other contributions and thus determines the sign of the total determinant as $n$ becomes large. \  It follows under the $k$-determinants assumption that $a_{1}a_{2} \cdots a_{k} > 0$. \  This completes the proof in the case in which the number of atoms is finite: asymptotically positive $k$-determinants force  $a_{1}a_{2} \cdots a_{k} > 0$.

Turning to the case in which the measure is countably atomic, we follow the above argument with improvements of Claims 4 and 5. \ Recall that we have assumed that $a = (a_i)_{i=1}^\infty$ is in $\ell^1$;  set $s := \sum_{i=1}^\infty |a_i|$.

\medskip

\noindent Claim 4$^\prime$.  For $\mathcal{C}$ with exactly one $j > k$ (say, $j_1$) and the remaining $k-1$ $j$ in $\{1, \ldots, k\}$ the contribution is bounded by the expression in \eqref{eq:boundonesmall}. \  Further, the sum of the contributions to the total determinant of all such $\mathcal{C}$ is bounded by
\begin{equation}  \label{eq:boundonesmallctbl}
s B^{k-1}r_1^n r_2^n \cdots r_{k-1}^n r_{k+1}^n k! r_1^{k((2k-2)}.
\end{equation}

\noindent Proof of Claim 4$^\prime$.  The first assertion is exactly as in Claim 4;  the second follows from summing these over all $j_1 > k$, using that $r_{k+1} \geq r_j$ for any $j \geq k+1$, and using the definition of $s$.

\medskip

\noindent Claim 5$^\prime$.  For $\mathcal{C}$ with at least two $j$ greater than $k$ -- say $j_1$ the smallest and $j_2$ the next smallest of these -- the contribution to the total determinant is bounded by the expression in \eqref{eq:boundfortwosmall}. \  For a fixed such $j_1$ -- the smallest of the $j$ larger than $k$ in $\mathcal{C}$ -- the sum of all such contributions for such $C$ is bounded by
\begin{equation}  \label{eq:boundspecjctbl}
|a_{j_1}| s B^{k-2}r_1^n r_2^n \cdots r_{k-2}^n r_{j_1}^n r_{j_2}^n k! r_1^{k(2k-2)}.
\end{equation}
Finally, the sum of all contributions to the total determinant from all $\mathcal{C}$ containing at least two $j$ larger than $k$ is bounded by
\begin{equation}  \label{eq:boundallctbl}
s^2 B^{k-2}r_1^n r_2^n \cdots r_{k-2}^n r_{j_1}^n r_{j_2}^n k! r_1^{k(2k-2)}.
\end{equation}

\noindent Proof of Claim 5$^\prime$.  The first assertion is exactly Claim 5;  the second assertion follows from fixing $j_1$, summing  the bounds in \eqref{eq:boundfortwosmall} over all $j_2$ larger than $j_1$ and using the definition of $s$. \  The third assertion follows by summing the bounds in \eqref{eq:boundspecjctbl} over all $j_1 > k$.

But then $D$ dominates, using \eqref{eq:boundonesmallctbl} and \eqref{eq:boundallctbl}, the sum of all other contributions to the total determinant and thus positivity of that determinant for $n$ large forces $D > 0$ and thereby $a_{1}a_{2} \cdots a_{k} > 0$. \  This completes the proof for the countably infinite case:  asymptotically positive $k$-determinants forces $a_{1}a_{2} \cdots a_{k} > 0$.

Finally, it is clear from the proof that the arguments are reversible to yield the converse claim. \pfend

\medskip
We now may obtain Theorem  \ref{th:CPDfirststatement} by turning our attention to $k$-hyponormality (in fact, in some asymptotic version). \  Observe that to make some definition like ``asymptotically positive $k$-hyponormality'' (requiring strict positivity of the moment matrices) is simply equivalent to the assumption that the shift $W_\alpha$ has some restriction to one of the canonical invariant subspaces $\bigvee \{e_i: i \geq I \}$ which is positively $k$-hyponormal, in the sense that the matrices are positive definite and not positive semi-definite. \  Note again that in our situation, we won't have $k$-hyponormality with determinants equal to zero, again by the rank arguments.

\begin{theorem}
Let $\mu$ be a (potentially signed) measure of the form $\sum_{i=1}^\infty a_i \delta_{r_i}$ where $a = (a_i)_{i=1}^\infty \in \ell^1$, $a_1 > 0$, and the $r_i$ form a strictly decreasing positive sequence. \  Suppose as well that $\sum_{i=1}^\infty a_i = 1$. \  Let $\gamma = (\gamma_n)_{n=0}^\infty$ be the moment sequence for $\mu$, so $\gamma_n = \int_{\mathbb{R}_+} t^n \, d \mu(t)$. \   Let $W_\alpha$ be the weighted shift with weight sequence $\alpha$ induced by the moments $\gamma_n$. \  Then the restriction of $W_\alpha$ to some canonical invariant subspace $\bigvee \{e_i: i \geq I \}$ is $k$-hyponormal for some $k=1, 2, $ if and only if the first $k+1$ non-zero $a_i$ are strictly positive.
\end{theorem}

\noindent Proof.  Since $k$-hyponormality forces positivity of the matrices of size less than $(k+1)\times(k+1)$ (and recall $a_1 > 0$), one deduces successively $a_1 a_2 > 0$, \ldots, $a_{1}a_{2} \cdots a_{k+1} > 0$, and the conclusion is obvious.  \pfend

Observe that, per the discussion of the construction of a counterexample as before Definition \ref{def:asymposkdets}, one may not deduce from the $a_i$ positive actual $k$-hyponormality;  what may fail is some $k$-hyponormality test for small $n$.

The reader may reasonably ask why we did not choose to take the ``restriction $k$-hyponormal'' as the useful condition as opposed to introducing ``asymptotically positive $k$-determinants.''  This question is rendered even more reasonable when one notes that a Berger-type charge $\mu$ for $W_\alpha$ and the Berger-type charge $\hat{\mu}$ for some canonical restriction are related by $\hat{\mu}(t) = t^n \mu(t)$ (for appropriate $n$, and with a normalization to a moment sequence beginning with $1$). \  The reason is that even in the absence of $k$-hyponormality, the signs of the determinants of the matrices considered yield information about the signs of the densities. \  Recall that $k$-hyponormality implies all the lesser $j$-hyponormalities, and therefore it is possible to have the determinants of the $k \times k$ Hankel moment matrices positive for a particular $k$ even without $k$-hyponormality. \  Thus we have the following as a corollary of the proof of the theorem -- we leave to the reader similar results -- and this provides motivation for the approach taken.

\begin{corollary}  \label{cor:densityinfonokHN}
Suppose $W_\alpha$ is a finite or countably atomic weighted shift with moment sequence having a Berger-type charge $\mu = \sum_{i=1}^\infty a_i \delta_{r_i}$ where $a = (a_i)_{i=1}^\infty \in \ell^1$, $a_1 > 0$, and the $r_i$ form a strictly decreasing positive sequence. \  Suppose as well that $\sum_{i=1}^\infty a_i = 1$, and assume for convenience that none of the $a_i$ is zero except possibly the tail (should the shift be finitely atomic). \  If $W_\alpha$ has both strictly positive $(k-1)$-determinants and $k$-determinants, then the density $a_k$ is positive. \  The same conclusion holds if $W_\alpha$ has both strictly negative $(k-1)$-determinants and $k$-determinants. \ If the determinants of the $(k-1) \times (k-1)$ and $k \times k$ Hankel moment matrices, for all $n$ large, have different (strictly positive or negative) signs, then the density $a_k$ is negative.
\end{corollary}

We remark that the question of how $k$-hyponormality affects non-atomic measures deserves study.

\subsection{Berger-type charges and CPD}   \label{BergerchargesandCPD}

The study of matrices which are conditionally positive definite is long-standing (see, for example, \cite{Bh}, \cite{HJ1}, and \cite{HJ2}). \    Recall that a matrix $M$ on $\mathbb{C}^n$ is positive definite (respectively, positive semi-definite) if for all $v = (v_1, \ldots, v_n)^T$ one has $v^T M v > 0$ (respectively, $v^T M v \geq 0$ for all $v$). \  The matrix $M$ is conditionally positive definite if $v^T M v > 0$ for all vectors $v$ satisfying the additional property $\sum_{i=1}^n v_i = 0$. \  (In practice, there is a slight abuse of language in which positive definite and conditionally positive definite include their ``semi-definite'' analogs.)

Conditionally positive definite matrices were used in a certain operator theory (weighted shift) context in \cite{BCE3} and there has been extensive study of conditionally positive definite (CPD) operators and weighted shifts in \cite{JJS} and \cite{JJLS}. \  We present next some results indicating the utility of Berger-type charges in this context.

Recall that in the study of the geometrically regular weighted shifts (GRWS) in \cite{BCE3} certain sectors of the Magic Square remained somewhat \textit{terra incognita}:  not $\mathcal{MID}$ or subnormal or even one of the $k$-hyponormalities and not completely hyperexpansive. \  We may make progress on understanding at least the additional subsector VIIIB by realizing that it consists of shifts which, although not in general CPD, can have their weights multiplied by a constant so as to become CPD. \  Recall from \cite[Theorem 2.2.13]{JJS} and surrounding discussion that a constant multiple of the weights of a CPD shift (there phrased in terms of moment sequences) need not result in a CPD shift. \ The next result identifies subsector VIIIB as those non-subnormal GRWS which have such a multiplier to CPD.

Subsector VIIIB is that subsector in the third quadrant of the Magic Square bounded above by the line $D = p N$ and below by $D = p^2 N$ (see Figure 1). \ It is known from the pattern of densities for points in this subsector that the induced shift is not normal. \  The next result identifies subsector VIIIB as those (and all of those) non-subnormal GRWS which have such a multiplier to CPD.

\begin{theorem}  \label{th:GRWSmultiptoCPD}
Let $p > 1$ and let $(N,D)$ in the Magic Square satisfy $-1 < p^2 N \leq D \leq p N < 0$, so as to yield a point in subsector VIIIB. \  We consider only multipliers such that after the multiplication and the application of $\Delta^2$, the new weights squared are non-negative. \ For such a pair $(N,D)$ not on the special line $D = p N$, there is exactly one non-zero constant multiple of the weights of $W_{\alpha(N,D)}$, namely by $\sqrt{p}$, so that the resulting shift is CPD. \  If the pair $(N,D)$ is on the special line in subsector VIIIB, there are exactly two non-zero multipliers, namely $1$ and $\sqrt{p}$, so that the resulting shift is subnormal. \  If $(N,D)$ is any point not in subsector VIIIB and not in Sectors I, II, and III, nor on the special lines in Sector IV (where the induced shifts are subnormal) there is no non-zero multiplier of the weights of $W_{\alpha(N,D)}$ to render it CPD.
\end{theorem}

\noindent Proof.  Let us consider the ``generic'' point in Subsector VIIIB on neither $D = p N$ nor $D = p^2 N$. \  The Berger-type charge for the induced shift $W_{\alpha(N,D)}$ is countably (not finitely) atomic, of the form
\begin{equation}  \label{eq:measureforS8}
C\left(\delta_1 + \sum_{i=1}^\infty c_i \delta_{\frac{1}{p^i}}\right),
\end{equation}
where $C$ is a normalizing constant, $c_1 = m_1 < 0$, $c_2 = m_1 m_2 > 0$, and the remaining $c_i$ are all positive. \  (See \cite{BCE3} for discussion of the $c_i$ and $m_i$, from which one computes the signs to yield the result in the list of density patterns on page \pageref{listofdensities}.)  Let $\gamma = (\gamma_n)_{n=0}^\infty$ be the moment sequence for $W_{\alpha(N,D)}$.

To multiply the weights by $\sqrt{k}$, for some $k > 0$, is to multiply the moments by $k^n$ so that the new moment sequence $\hat{\gamma}$ satisfies $\hat{\gamma}_n =  k^n \gamma_n$ for all $n$. \  This may be obtained by convolution with the measure $\delta_k$, using that the convolution of atomic measures $\delta_a$ and $\delta_b$ is $\delta_{a b}$. \  It is easy to compute that to obtain the measure for $\Delta^2$ of a single atom $\delta_a$ one should use $(1-a)^2 \delta_a$. \  Multiplying, taking $\Delta^2$, and summing over the atoms of \eqref{eq:measureforS8}, we obtain that the result of the multiplication of the weights by $\sqrt{k}$ and subsequent $\Delta^2$ is
\begin{equation}  \label{eq:measurepostkDeltaup2}
C\left((k-1)^2\delta_k + \sum_{i=1}^\infty c_i (1 -\frac{k}{p^i})^2 \delta_{\frac{k}{p^i}}\right).
\end{equation}
Since we want the measure for a subnormal shift, we must choose $k$ so that the density at $\frac{1}{p}$, which is $c_1 (1 - \frac{k}{p})^2$, becomes zero else it will be negative. \  This forces $k = p$, and clearly no other $k$ works.

For a point on the line $D = p^2 N$ the computation is the same, although in this case the original measure is three-atomic of the form
$$C\left(\delta_1 + c_1 \delta_{\frac{1}{p}} + c_2 \delta_{\frac{1}{p^2}}\right).$$

Consider now a point on the line $D = p N$, and we will perform an alternate computation to illustrate it and to prepare for some subsequent discussion. \  One may show that to achieve the effect of $\Delta$ on a \textbf{weight} sequence $\alpha$, what results is weights of the form
\begin{equation}   \label{eq:weighttransformDelta}
\sqrt{\alpha_n^2\frac{(\alpha_{n+1}^2 - 1)}{(\alpha_{n}^2 - 1)}},\quad n = 0, 1, \ldots.
\end{equation}
One computes via two applications of this formula, following a multiplication by $1$, that the resulting weights are constant at $\sqrt{\frac{1}{p}}$ obviously yielding a subnormal shift;  alternatively, after multiplication by $\sqrt{p}$ and $\Delta^2$ the resulting weights are
$$\sqrt{p}\sqrt{\frac{(p^2)^n + \frac{M^2}{p}}{(p^2)^n + M^2 p}},$$
and this is a multiple of a shift from the special line $D = p^2 N$ with parameter $p^2$ and such shifts are subnormal (see \cite{BCE3}). \  Returning to the measure approach, one may compute as before that no other multiplier will work.

Finally, for points outside of Subsector VIIIB and the subnormal Sectors I, II, and III, and the special lines in Sector IV, the measure calculation of multiplication followed by $\Delta^2$ leaves one with multiple negative densities (along with some positive densities). \  These negative densities cannot all be made zero by any single choice of $k$, so the resulting shift is not subnormal so even a multiplication of the original cannot be CPD.  \pfend

\begin{remark}  \label{remark:signedrepresentations}
There is a crucial point to be made in connection with the proof just given;  further examination of the proof will underline the need to accord recognition to signed representing measures (Berger-type charges). \  Suppose that, for a point on the line $D = pN$, one were to do the measure calculation as in the first part of the proof and in the (unmultiplied) case $k = 1$. \  On this line, the measure is two-atomic of the form
$$
C(\delta_1 + m_1 \delta_{\frac{1}{p}})
$$
where $m_1 < 0$. \ Taking $\Delta^2$ yields the measure
$$C\cdot m_1 \cdot (1 - \frac{1}{p})^2 \delta_{\frac{1}{p}}.$$
Here is the point:  although the density is negative, one cannot discard this as failing to yield a subnormal shift. \  The key is that if we can be brave, and produce weights by ratios of successive moments, then as far as weights go this (negative) measure produces perfectly good subnormal weights. \  In fact, these weights are exactly those produced by the (positive) measure
$$-C\cdot m_1 \cdot (1 - \frac{1}{p})^2 \delta_{\frac{1}{p}}.$$

Once stated this is obvious:  if $\mu$ is any Berger measure, computation of weights squared from $-\mu$ by ratios of moments yields the same weights and thus subnormality. \  Traditionally a moment sequence $\gamma$ for weighted shifts had $\gamma_0  = 1$;  in practice, as one can note in the proof of Theorem \ref{th:GRWSmultiptoCPD} in which we examine a measure not producing a moment sequence with $\gamma_0  = 1$, we sometimes let this pass by, perhaps muttering ``well, Stieltjes sequence.''  But we must go farther:  the only way to resolve the apparent contradiction of the previous computations is to accept the negative measure $C\cdot m_1 \cdot (1 - \frac{1}{p})^2 \delta_{\frac{1}{p}}$ as being a satisfactory representing measure for a subnormal shift.

It follows as a second point that we must distinguish carefully between whether we are considering moment sequences (in which case surely $\gamma \neq -\gamma$) or discussing weighted shifts, in which the weights are fundamental and the moment sequences auxiliary. \  We will have more to say about this shortly.

When do negative densities indicate a failure of subnormality?  Answer:  if we have a Berger-type charge $\mu$ with both (strictly) positive and (strictly) negative densities, then the resulting shift is not subnormal. \  But we are compelled to accept a Berger-type charge with all negative densities as being a perfectly good (signed) representing measure for a subnormal shift.

The careful reader may check that elsewhere in this paper when we have taken a negative density to indicate failure of subnormality, there have always been the needed positive densities to render the conclusion correct.
\end{remark}

To aid in the distinction between the study of moment sequences and the study of shifts, we will say that a shift is ``CPD-weights'' if $\Delta^2$ of the moment sequence yields the moment sequence for a subnormal shift, reserving ``CPD'' for when $\Delta^2$ yields a positive definite moment sequence.

There is a generalization of the previous result from the GRWS situation to more general countably atomic Berger-type charges. \  The proof of the following -- a more careful statement of Theorem \ref{th:CPDfirststatement} -- is a modification of the argument for the proof of Theorem \ref{th:GRWSmultiptoCPD}, and we omit it.

\begin{theorem}
Let $\mu$ be a countably atomic (potentially signed) measure with support in $[0, \infty)$. \  If $\mu$ has all negative or all positive densities, the shift $W_\alpha$ is subnormal and hence conditionally positive definite in the sense of CPD-weights. \  If $\mu$ has two atoms with one density of each sign, then there are exactly two (non-zero) constants which, multiplying the weights of $W_{\alpha}$, make it CPD-weights. \  If $\mu$ has more than two atoms, with either all positive densities but one or all negative densities but one, there exists exactly one (non-zero) multiplier of the weights of $W_\alpha$ producing a CPD-weights shift. \  If $\mu$ has at least two atoms with positive densities and at least two with negative densities, no non-zero multiple of $W_\alpha$ is CPD-weights.
\end{theorem}

We close with a discussion of completely hyperexpansive shifts:  in what sense, if any, are these CPD?  It is well known that if $\gamma$ is the (traditionally, positive) moment sequence of a CHE shift, then $\Delta(\gamma)$ is completely monotone (see \cite[Ch. 4, Lemma 6.3]{BCR}). \ If we take $\nabla$ of a completely monotone sequence, it is again completely monotone, but $\Delta$ of a completely monotone sequence is the negative of a completely monotone sequence. \  Thus a CHE shift (except in trivial cases such as the unweighted unilateral shift) is not CPD. \  However, it is CPD-weights, since from the point of view of weights there is no difference between the sequence of operations $\Delta$ followed by $\nabla$ as opposed to $\Delta^2$.

We also point out that there is, with the aid of a signed measure, a ``CPD-like'' representation for the moment sequence of some CHE shifts. \  Recall that the (traditional, positive) moment sequence of a CHE shift has the L\'{e}vy-Khinchin representation of the form
$$
\gamma_n = 1 + b n + \int_{[0,1)} (1-x^n) \, d \nu(x) \quad (n = 0, 1, 2, \ldots),
$$
where $b \in \mathbb{R}$ and $\nu$ is a positive Radon measure on $[0,1)$ (\cite[Ch. 4, Prop. 6.11]{BCR}). \ Equivalently, 
$$
\gamma_n = 1 + \mu(\{1\})\cdot n + \int_{[0,1)} \frac{1-x^n}{1-x} \, d \mu(x) \quad (n = 0, 1, 2, \ldots),
$$
where $\mu$ is the finite (positive) Borel measure on $[0,1]$ associated with $\Delta\gamma$ (see \cite[page 3747]{At}. \ For each $n = 0, 1, 2, \ldots$, we write
\begin{eqnarray*}
\gamma_n &=& 1 + \mu(\{1\})\cdot n + \int_{[0,1)} \frac{1-x^n}{1-x} \, d \mu(x) \\
&=& 1 + b n +  \int_{[0,1)} \frac{x^n - 1}{x-1} \, d \mu(x) \\
&=& 1 + b n + \int_{[0,1)} \frac{x^n - 1}{(x-1)^2} \cdot (x-1) d \mu(x) \\
&=& 1 + (b + \nu([0,1))) \cdot n +  \int_{[0,1)} \frac{x^n - 1}{(x-1)^2} - \frac{n(x-1)}{(x-1)^2} \cdot (x-1) d \mu(x)\\
&=& 1 + (b + \nu([0,1))) \cdot n +  \int_{[0,1)} \frac{x^n - 1 - n(x-1)}{(x-1)^2} \cdot (x-1) d \mu(x)\\
&=& 1 + (b + \nu([0,1))) \cdot n +  \int_{[0,1)} Q_n(x) \cdot (x-1) d \mu(x),
\end{eqnarray*}
where $Q_n(x) := \dfrac{x^n - 1 - n(x-1)}{(x-1)^2}$ is as in the CPD representation in \cite{JJS}. \ The expression above gives a ``CPD-like'' representation, at the cost of a negative representing measure $(x-1) \mu(x)$. \ From \cite[Lemma 2.2.1]{JJS} we have that $\Delta^2(Q(\cdot))(x)_n = x^n$ for $n \in \mathbb{Z}_+$ and $x \in \mathbb{R}$, we can see how a negative Berger-type charge for the moment sequence of a subnormal shift will arise naturally from $\Delta^2$ applied to the above representation. \  This both indicates potential new objects of study and emphasizes the distinction between CPD and CPD-weights.

\medskip

\noindent \textbf{Acknowledgments.} \ The authors wish to express their appreciation for support and warm hospitality during various visits (which materially aided this work) to Bucknell University, the University of Iowa, and the Universit\'{e} des Sciences et Technologies de Lille, and particularly the Mathematics Departments of these institutions. \ The second named author was partially supported by NSF grant DMS-2247167. \ Several examples in this paper were obtained using calculations with the software tool \textit{Mathematica} \cite{Wol}.

\end{document}